\documentclass[10pt]{article}
\usepackage{amsthm,amsmath,amsfonts}
\usepackage{fullpage,enumerate}
\usepackage{tikz}
\usepackage[colorlinks,citecolor=blue,urlcolor=blue]{hyperref}
\usepackage{hypernat}
\usepackage[notref,notcite]{showkeys}

\newcommand{\la}{\lambda}

\begin{document}

\numberwithin{equation}{section}

\title{Does random dispersion help survival?}

\maketitle

Rinaldo B. Schinazi

Department of Mathematics

University of Colorado, Colorado Springs

rschinaz@uccs.edu

\bigskip

{\bf Abstract}  Many species live in colonies that prosper for a while and then collapse. After the collapse the colony survivors disperse randomly and found new colonies that may or may not make it depending on the new environment they find.  We use birth and death chains in random environments to model such a population and to argue that random dispersion is a superior strategy for survival.

\section{ The model}

Many species (such as ants or early humans) survive according to the following pattern. A few individuals start a colony, the colony prospers until it collapses due to exhaustion of resources or some external factor, a few individuals survive the collapse, disperse randomly and found new colonies. The pattern starts over. In this paper we are interested in modeling this process and comparing it to other survival strategies.

Next we give an informal description of our stochastic model. We will be more precise in the proof of Theorem 1.

The individuals are divided in independent colonies. Each colony is a birth and death chain with birth rate $\la$ sampled from a fixed distribution $\mu$ and death rate 1. Moreover, each colony is associated with a random time $T$ sampled from a fixed distribution $\nu$. At time $T$ each individual of the colony starts its own independent colony. Each new colony has a new $\la$ and a new $T$ sampled from $\mu$ and $\nu$, respectively. The new $(\la,T)$ is independent of everything else but $\la$ and $T$ need not be independent.

We are particularly interested in distributions $\mu$ for which $\mu([0,1))$ is almost 1. That is, most of the $\la$'s that are sampled are below 1 and the new colony collapses very quickly (this also models the fact that most individuals die when their colony collapses).  Can such a population survive? We will show that it can under certain conditions.

\medskip

{\bf Theorem 1. } {\sl Consider a process in which the individuals are divided in a random number of independent colonies. Each colony is a birth and death chain for which the birth rate is sampled from a fixed distribution $\mu$ and the death rate is 1. Each colony is associated with a random time $T$ which is sampled from a fixed distribution $\nu$. At time $T$ each individual in the colony starts its new independent colony with a new $\la$ and a new $T$. This population survives if and only if
$$E[\exp((\Lambda-1)\tau)]>1,$$
where $\Lambda$ has distribution $\mu$ and $\tau$ distribution $\nu$.}

\medskip

As we show next an easy consequence of Theorem 1 is that if $E(\Lambda)\geq 1$ then the population survives.

\medskip

{\bf Corollary 1.} {\sl Assume that $\Lambda$ and $\tau$ are independent. If $E(\Lambda)>1$ then the population survives. Moreover, if $E(\Lambda)=1$ and $\Lambda$ is not constant a.s. then the population also survives.}

\medskip

Comparing our model to a birth and death chain with fixed birth rate $E(\Lambda)$ and death rate 1 we see that if the fixed birth rate process survives so does our process. Moreover, we know that a critical birth and death chain dies out (see Schinazi (1999) or Karlin and Taylor (1975)). On the other hand if $E(\Lambda)=1$ then our model survives provided the birth rates are sampled from a non degenerate distribution. Hence, random dispersion seems to help survival. We will actually give an example for which $E(\Lambda)$ can be arbitrarily close to 0 and our model still survives. 

\medskip

{\it Proof of Corollary 1.} 

By Jensen's inequality we have
$$E(\exp(\Lambda-1)\tau)\geq \exp[E((\Lambda-1)\tau)].$$
There are two cases to consider.

$\bullet$ If $E(\Lambda)>1$ then the r.h.s. is strictly bigger than 1 (provided $E(\tau)>0$) and therefore the process survives.

$\bullet$ Assume $E(\Lambda)=1$. It is easy to check that Jensen's inequality is strict provided $(\Lambda-1)\tau$ is not a constant. Hence, assuming $\Lambda\not= 1$ a.s. we have a strict inequality and therefore $E(\exp(\Lambda-1)\tau)>1$. Thus, the process survives.

\subsection{An example}

Consider the case where $\mu$ concentrates on only two values $\lambda_1$ and $\lambda_2$ with probability $p$ and $1-p$, respectively. We also assume that $\nu$ is an exponential distribution with parameter $a$ and that $\Lambda$ and $\tau$ are sampled independently. Let $(Y_t)_{t\geq 0}$ be the total number of individuals at time $t$.

There are several cases to consider.

$\bullet$ Assume that $0\leq \la_1<\la_2\leq 1$. We can couple $(Y_t)_{t\geq 0}$ with a classical birth and death chain $(W_t)_{t\geq 0}$ with birth rate $\la_2$ and death rate 1 in such a way that $Y_0=W_0=1$ and $Y_t\leq W_t$ for all $t\geq 0$. This is so because the death rates are the same but the birth rate for $(W_t)_{t\geq 0}$ is larger or equal than for $(Y_t)_{t\geq 0}$. Since we assume $\la_2\leq 1$ we know that $(W_t)_{t\geq 0}$ dies out with probability 1. Therefore, so does $(Y_t)_{t\geq 0}$.

$\bullet$ Assume that $1<\la_1<\la_2$. We can couple the model in random environment with a classical birth and death chain with birth rate $\la_1$. In this case the model in random environment dominates the classical model. The latter model survives with positive probability since $\la_1>1$. Hence, the model in random environment survives as well.

$\bullet$ The interesting case is $0\leq \la_1\leq 1<\la_2$ because now one environment is subcritical and the other is supercritical. We turn to this case for the rest of this section.

Let 
$$m=E(\exp(\Lambda-1)\tau)=pE(\exp(\la_1-1)\tau)+(1-p)E(\exp(\la_2-1)\tau)$$
where $\tau$ has an exponential distribution with parameter $a$.

If $\la_2\geq a+1$ then $m=+\infty$ and $(Y_t)_{t\geq 0}$ survives with positive probability.  Since we are free to pick any $\la_1\geq 0$ and any $p$ in $(0,1)$ we see that $E(\Lambda)$ can be arbitrarily close to 0 and the process still survives.

Assume now that $\lambda_2<a+1$.
We have
$$m=p\frac{a}{a+1-\lambda_1}+(1-p)\frac{a}{a+1-\lambda_2}.$$
A little algebra shows that $m>1$ if and only if
$$a(-1+p\la_1+(1-p)\la_2)>(1-\la_1)(1-\la_2).$$
Note that 
$$E(\Lambda)=p\la_1+(1-p)\la_2.$$
There are three cases to consider.

$\bullet$ If $E(\Lambda)>1$ then $m>1$ if and only if
$$a>\frac{(1-\la_1)(1-\la_2)}{-1+E(\Lambda)}.$$
Since the r.h.s. is negative this is true for all $a> \la_2-1$. We also know that $(Y_t)_{t\geq 0}$ survives for $a\leq \la_2-1$. Hence, if $E(\Lambda)>1$ the process survives for all $a>0$. 

$\bullet$ If $E(\Lambda)<1$ then $m>1$ if and only if
$$a<a_c=\frac{(1-\la_1)(1-\la_2)}{-1+E(\Lambda)}.$$
Hence, $(Y_t)_{t\geq 0}$ survives for $\la_2-1<a<a_c$. We also have survival for $a\leq \la_2-1$. Hence, if $E(\Lambda)<1$ the process survives if and only if $a<a_c$.

$\bullet$ If $E(\Lambda)=1$ then  $m>1$ if and only if 
$$0>(1-\la_1)(1-\la_2).$$
This is true for all $a>0$. Hence, the process survives. Note that the process does not survive if $\la_1=\la_2=1$. That is, the critical process survives if and only if the environment is truly random.

\medskip

The results for $E(\Lambda)\geq 1$ were already known by Corollary 1. The case $E(\Lambda)<1$ is more interesting. What determines survival in this regime is the parameter $a$ and hence the distribution of $\tau$.  Interestingly, the process survives provided the changes in environment do not happen too frequently (i.e. $a$ has to be less than $a_c$) even though one of the two environments can be quite hostile (such as $\la_1=0$).

\section{ A null model}

The preceding example shows that random dispersion helps survival as compared to a model with fixed environment. However, it is not clear whether easier survival is due to the multiplication of random environments or simply to a changing random environment. In this section we introduce a model for which there is no dispersion. The environment changes globally for the whole population. We will show that a globally changing environment does not help survival. It is in fact the combination of dispersion and changing environments that helps survival.

Let $\tau_1, \tau_2,\dots $ be a sequence of i.i.d. random variables. Define for $k\geq 1$
$$T_k=\tau_1+\tau_2+\dots+\tau_k.$$
Let $\mu$ be a probability distribution whose support is in $[0,+\infty)$.
Let $(X_t)_{t\geq 0}$ be a continuous time chain with the following rules. Set $X_0=1$ and let $\lambda_0>0$ be sampled from $\mu$. For $t<T_1$ individuals give birth at rate $\lambda_0$ and die at rate 1. At time $T_1$, $\lambda_1$ is sampled from $\mu$ independently of everything else and is kept until time $T_2$. Between times $T_1$ and $T_2$ individuals give birth at rate $\lambda_1$ and die at rate 1. At time $T_2$ a new $\lambda_2$ is sampled and so on. In other words $(X_t)_{t\geq 0}$ is a birth and death chain with death rate fixed at 1 and birth rates changing at times $T_k$ for $k\geq 1$. Next we give a necessary and sufficient condition for survival.

\medskip

{\bf Theorem 2.} {\sl Consider a birth and death chain in random environment for which the death rate is 1 and the birth rate changes at times $T_k=\tau_1+\tau_2+\dots+\tau_k$ where $\tau_1, \tau_2,\dots $ are i.i.d. Birth rates are sampled from a fixed distribution $\mu$. Assume that $E(\tau_1)<+\infty$ and $E(\Lambda)<+\infty$ where $\Lambda$ has distribution $\mu$. Assume that the birth rates and the times $\tau_i$ are independent. Then, the process survives if and only if $$E(\Lambda)>1.$$
}

\medskip

Theorem 2 proves our point. The birth and death chain in a global random environment survives if and only if the classical birth and death chain with birth rate $E(\Lambda)$ survives. 
That is, without dispersion the random environment does not help survival. 

\medskip

\section{Proof of Theorem 1}

We start the process $(Y_t)_{t\geq 0}$ with a single individual. This first individual is our generation 0. We sample a random time $\tau_{0,1}$ from distribution $\nu$. The process $(Y_t)_{t\geq 0}$ for $t<\tau_{0,1}$ is a birth and death chain with birth rate $\lambda$ sampled from a distribution $\mu$ and a death rate equal to 1. The $j\geq 0$ individuals present at time $\tau_{0,1}$ constitute generation 1.   Assuming $j\geq 1$, each of these $j$ individual samples a new $\lambda$ from $\mu$ independently  of each other and starts its own colony. We associate to each of the $j$ colonies an exponential random variable $\tau_{1,k}$ for $k=1,\dots,j$.  The individuals present at times $\tau_{1,k}$ for $k=1,\dots,j$ constitute generation 2.  Each generation 2 individual samples a new $\lambda$, a new $\tau$  and starts a new colony and so on.

In summary, we have independent colonies attached to random times $\tau_{i,j}$. Each colony evolves as a birth and death chain with a fixed birth rate and a death rate equal to 1. At time $\tau_{i,j}$ each individual present in the colony starts its own independent colony and so on.  Each colony has a fixed birth rate which is sampled independently of everything else from the same  distribution $\mu$. The times $(\tau_{i,j})_{i\geq 0, j\geq 1}$ are i.i.d.

We will find a necessary and sufficient condition for survival.
We first define an auxiliary discrete time process $(V_n)_{n\geq 0}$. The random variable $V_n$ counts the number of individuals in generation $n$.
Hence, $V_0=1$. Let $n\geq 1$, assume that $V_n=j\geq 1$. Each of these $j$ individuals will start a new colony (at different random times). To each colony is attached an exponential time $\tau_{n,k}$, $k=1,\dots,j$.  Let $V_{n+1,k}$ be the number of individuals in colony $k$  present at time $\tau_{n,k}$.  We think of $V_{n+1,k}$ as being the offspring of the individual who started the $k$-th colony. Define
$$V_{n+1}=\sum_{k=1}^j V_{n+1,k}.$$
If for some $n\geq 0$ we have $V_n=0$ we set $V_k=0$ for all $k>n$. 

We claim that $(V_n)_{n\geq 0}$ is a (classical) Galton-Watson process. This is so because 
each individual in generation $n$ has the same offspring distribution as the individual in generation 0. Moreover, different individuals have independent offspring distributions. That is, the random variables $V_{n,k}$ are i.i.d. 
Hence, the process $(V_n)_{n\geq 0}$ survives if and only if
$$E(V_1)>1.$$
Given $\Lambda=\lambda$ and $\tau_{0,1}=t$, $V_1$ is the number of individuals at time $t$ of a birth and death process with birth rate $\lambda$ and death rate 1. Therefore,
$$E(V_1|(\la,t))=e^{(\lambda-1)t}.$$
Hence,
$$E(V_1)=E[\exp((\Lambda-1)\tau)].$$
Note that $(Y_t)_{t\geq 0}$ survives if and only if $(V_n)_{n\geq 0}$ survives.
This completes the proof of Theorem 1.

\section{ Proof of Theorem 2}

Recall that the number of individuals at time $t$ for this model is denoted by $X_t$. We define a discrete time process $(Z_n)_{n\geq 0}$ in the following way. Let $Z_0=1$ and for $k\geq 1$

$$Z_k=X_{T_k}.$$

We claim that the process  $(Z_n)_{n\geq 0}$ is a (discrete time) branching process in random environment, see Smith and Wilkinson (1969).  To see this consider the individuals present at time $T_k$. A new birth rate $\la$ and a new random time $\tau_{k+1}$ are sampled at time $T_k$. The pair $(\lambda,\tau_{k+1})$ determines the random environment between time $T_k$ and $T_{k+1}$.  Given a random environment the progenies at time $T_{k+1}$ of each individual present at time $T_k$ are independent and identically distributed.  This fits the definition of Smith and Wilkinson (1969).  In order to state their theorem we need more notation. 

The process $(Z_n)_{n\geq 0}$ starts with a single individual (i.e. $Z_0=1$). Assume that the first random environment is $(\tau,\Lambda)=(t,\la)$. Then $Z_1$ is the number of individuals  in a birth and death chain $B(t,\la)$ with birth rate $\lambda$ and death rate 1.

Let
$$b(t,\la)=E[B(t,\la)].$$
 We know that
$$b(t,\la)=e^{(\la-1)t}.$$

Smith and Wilkinson (1969) proved that under the integrability conditions 
$$E|\ln b(\tau,\Lambda)|<+\infty\mbox{ and }E|\ln P(B(\tau,\Lambda)\geq 1)|<+\infty$$
the necessary and sufficient condition for survival is given by
$$E(\ln b(\tau,\Lambda))>0.$$
We now check the integrability conditions. We have
$$E|\ln b(\tau,\Lambda)|=E|(\Lambda-1)\tau|.$$
This expected value is finite since $\Lambda$ and $\tau$ are independent and their first moments are finite. 

We now turn to the second integrability condition. One way for a birth and death chain not to be extinct by time $t$ is for the initial individual not to have died by time $t$. Hence,
$$P(B(t,\la)\geq 1)\geq e^{-t}.$$
Therefore,
$$|\ln P(B(t,\la)\geq 1)|\leq t.$$
Hence,
$$E|\ln P(B(\tau,\Lambda)\geq 1)|\leq E(\tau)<+\infty.$$

We are now ready to apply the necessary and sufficient condition for survival to the process $(Z_n)_{n\geq 0}$. We have
$$E(\ln b(\tau,\Lambda))=E[(\Lambda-1)\tau]=E(\Lambda-1)E(\tau).$$
We see that $(Z_n)_{n\geq 0}$ survives if and only if $E(\Lambda)>1$. Note that the process $(X_t)_{t\geq 0}$ survives if and only if the process $(Z_n)_{n\geq 0}$ survives. This concludes the proof of Theorem 2.

\bigbreak

{\bf References}

S. Karlin and Taylor (1975) {\it A first course in stochastic processes}, second edition.
Academic Press, New York.

R.B. Schinazi (1999) {\it Classical and spatial stochastic processes}, Birkhauser.

W.L. Smith and W.E. Wilkinson (1969). On branching processes in random environments. The Annals of Mathematical Statistics 40, 814-827.

\end{document}